\documentclass{amsart}

\usepackage{amssymb}

\newtheorem{thm}{Theorem}%[section]

\newtheorem{lem}[thm]{Lemma}

\newtheorem{prop}[thm]{Proposition}

\newtheorem{conj}[thm]{Conjecture}
\theoremstyle{definition}
\newtheorem{defn}[thm]{Definition}

\newtheorem{say}[thm]{}
\newtheorem{exmp}[thm]{Example}

\newtheorem{prob}[thm]{Problem}
\newtheorem{rem}[thm]{Remark}          

\newtheorem{ack}{Acknowledgments}

\newtheorem{defn-thm}[thm]{Definition--Theorem}  %!!!!!!!!!!!!!!!!!!!!!!!!
\newtheorem{defn-lem}[thm]{Definition--Lemma}  %!!!!!!!!!!!!!!!!!!!!!!!!
\theoremstyle{remark}

\renewcommand{\c}[0]{{\mathbb C}}  

\renewcommand{\o}[0]{{\mathcal O}} 
\newcommand{\z}[0]{{\mathbb Z}}

\renewcommand{\a}[0]{{\mathbb A}}

\newcommand{\p}[0]{{\mathbb P}}
\newcommand{\f}[0]{{\mathbb F}}
\newcommand{\q}[0]{{\mathbb Q}}

\newcommand{\qtq}[1]{\quad\mbox{#1}\quad}
\newcommand{\spec}[0]{\operatorname{Spec}}

\newcommand{\inter}[0]{\operatorname{Int}}    
\newcommand{\sing}[0]{\operatorname{Sing}}

\newcommand{\stab}[0]{\operatorname{stab}}

\newcommand{\rdown}[1]{\lfloor{#1}\rfloor}

\newcommand{\lcm}[0]{\operatorname{lcm}}

\def\into{\DOTSB\lhook\joinrel\to}

\begin{document}
\bibliographystyle{amsalpha}

\title[Bogomolov--Miyaoka--Yau inequality]{Is there a topological\\
Bogomolov--Miyaoka--Yau inequality?}
\author{J\'anos Koll\'ar}

\today

\maketitle

%\tableofcontents

Let $S$ be a smooth, complex,  projective, minimal  surface
of general type. The Bogomolov--Miyaoka--Yau inequality states
that $c_1(S)^2\leq 3 c_2(S)$ \cite{bog-4, rei-bog, miy-3, yau}.
In this note I want to address the
following question:
$$
\mbox{\it Is there a  topological analog of the
Bogomolov--Miyaoka--Yau inequality?}
$$
\noindent The $11/8$-conjecture 
\cite{matsu, furuta} can be viewed as such, but 
in Section 1 I write down another
possible variant. 
 Section 2 explores its relationship with the
Montgomery--Yang problem on  differentiable circle actions on
$S^5$ and Section 3 examines its connection
with the  H--cobordism group of 3--manifolds.
Related examples  and questions on algebraic surfaces 
are discussed in Section 4.
The last section studies the remarkable series of
hypersurfaces
$$
(x_1^{a_1}x_2+x_2^{a_2}x_3+ \cdots + x_{n-1}^{a_{n-1}}x_n+ x_n^{a_n}x_1=0)
$$
in weighted projective spaces. 

The aim of this note is two fold.
On the one hand, I would like to call attention to
several questions about algebraic surfaces with quotient
singularities that have interesting connections with
the topology of 3-- or 5--manifolds.
The study of algebraic orbifolds led to many interesting
examples in  topology and differential geometry 
(see, for instance, \cite{briesk, or-wa} or the 
 recent papers \cite{bgk, e-surf}),
 but there should be many
more connections.

On the other hand, more speculatively, I hope that
methods developed around the  Bogomolov--Miyaoka--Yau inequality
can be adapted to the topological setting,  leading to
progress on the questions mentioned in Sections 2 and 3.

\section{The orbifold Bogomolov--Miyaoka--Yau inequality}

The Bogomolov--Miyaoka--Yau inequality can be
generalized to orbifolds. These are normal projective surfaces
whose singularities are locally analytically isomorphic to
 quotient singularities $\c^2/G$ where $G\subset GL(2,\c)$
is a finite group whose action is fixed point free outside the origin.
Such a surface has finitely many singular points,
and locally at each of them $S$ is topologically the cone
over a 3--manifold $S^3/G$ where $G\subset U(2,\c)$
is a subgroup acting without fixed points.
This 3--manifold is called the
{\it link} of $s\in S$ and it is denoted by $L_s$.
(The version where a surface is allowed to have orbifold
structure in codimension one is also interesting \cite{bgk, e-surf}, but
it will not be considered here.)

An easy way to get such an orbifold is to take the quotient of  a smooth 
projective surface $X$ by a finite group $G$ acting on $X$
with only isolated fixed points.
The most interesting examples are, however, those that do not
arise as a global quotient. 
  For a complex surface, $c_2(S)$ is
the same as the topological Euler characteristic.
Following Thurston,
we introduce the
{\it orbifold Euler characteristic}
$$
e_{orb}(S):=e(S)-\sum_{s\in \sing S} \Bigl( 1-\frac1{|\pi_1(L_s)|}\Bigr).
$$
The orbifold version of the
Bogomolov--Miyaoka--Yau inequality is the following,
developed in the series of papers
\cite{sak, miy, kns, megy}. 

\begin{thm}\label{obmy.thm}  Let $S$ be a normal projective surface with
quotient singularities such that $-c_1(S)$
is ample (or at least nef). Then
$$
c_1(S)^2\leq 3e_{orb}(S).
\eqno{(\ref{obmy.thm}.1)}
$$
\end{thm}

Since  $-c_1(S)$
is  nef,  $c_1(S)^2\geq 0$, thus we also get
the following weaker inequality,
which also holds when $c_1(S)$ is nef by \cite{ke-mc}
$$
0\leq e_{orb}(S).
\eqno{(\ref{obmy.thm}.2)}
$$
While (\ref{obmy.thm}.2) gives nothing interesting 
for smooth surfaces, it has very interesting
consequences for singular surfaces.

 Let $S$ be a normal projective surface with
quotient singularities whose 
 second Betti number  is 1. Then either $\pm c_1(S)$
is ample or $c_1(S)=0$. Thus (\ref{obmy.thm}.2) applies
and so $0\leq e_{orb}(S)$. 
It is easy to see (by looking at the Albanese map)
that $b_1(S)=0$, thus $e(S)=3$ and so $S$ is a
rational homology $\c\p^2$. Thus we get to following 
 expanded version
of  (\ref{obmy.thm}.2):
$$
\sum_{s\in \sing S} \Bigl( 1-\frac1{|\pi_1(L_s)|}\Bigr)\leq  3.
\eqno{(\ref{obmy.thm}.3)}
$$
In particular, $S$ has at most 6 singular points.
\medskip

Let us now make the bold (or foolish)
guess that the inequality
(\ref{obmy.thm}.3) is topological in nature.
It may be more convenient to formulate the conjecture
for smooth, compact 4--manifolds $M$ whose
boundary components are spherical, that is, their universal cover is $S^3$.
One can then attach cones to each boundary component
to get a 4--dimensional orbifold $S$.
We are mainly interested in the cases when $S$ is a
 homology $\c\p^2$. Correspondingly,
$H_1(M,\z)=0$ and $H_2(M,\z)\cong \z$.

\begin{conj}[Smooth Bogomolov--Miyaoka--Yau inequality]
\label{tbmy.thm}
Let $M^4$ be a smooth, compact 4--manifold with spherical
boundary components $\partial M^4=\cup_i L_i$. Assume that
$H_1(M,\z)=0$ and $H_2(M,\z)\cong \z$. Then
$$
\sum_i \Bigl( 1-\frac1{|\pi_1(L_i)|}\Bigr)< 3.
\eqno{(\ref{tbmy.thm}.1)}
$$
In particular, $M$ has at most 5 boundary components.
\end{conj}

\begin{say}[The origins of the conjecture] While, to the best of my knowledge,
the above conjecture is new, various forms of it have been implicit
in the works of several authors.

Controlling the singularities of a variety by global invariants has
been one of the  aims of the papers
\cite{miy,  megy, ke-mc}. It is quite natural to go from
an algebraic formulation to a purely topological one as above.

On the topological side, the Montgomery--Yang problem on circle actions
on 5--manifold was explicitly studied as a problem  on 4--dimensional
orbifolds in \cite{fi-st}. As we discuss in Section 2, the
relationship with (\ref{tbmy.thm}) is close but the two problems are
not identical.

The history of the topological work on 
H--cobordisms of Seifert fibered 3--manifolds
is reviewed in \cite{sav}.
The same topic appears in singularity theory as the study of smoothings
of quasihomogeneous surface singularities. The link of such a singularity
is a Seifert fibered 3--manifold, and each smoothing exhibits the link as
the boundary of the Milnor fiber. 
In this setting, (\ref{tbmy.thm}) is closely related to 
\cite[3.10]{wahl-my}. His results imply (\ref{tbmy.thm})
for smoothings of negative weight.
\end{say}

\begin{say}[Comments and possible generalizations] 

1. In the algebraic case,
if equality holds in (\ref{obmy.thm}.3) then $c_1(S)=0$ and
$S$ is the quotient of an Abelian or a K3 surface by a cyclic group.
In particular, $H_1(S\setminus \sing S,\z)\neq 0$. This is why
the $\leq$ in (\ref{obmy.thm}.3)
 was changed to a strict inequality in (\ref{tbmy.thm}.1).
An algebraic rational homology $\c\p^2$   can never have 6 singular points
and I don't know any examples with 5 singular points.
There are, however, many examples with 4 singular points
(\ref{4pt.icosa}).

2. Unlike in the algebraic case, the restriction
$H_1(M,\z)=0$ is not a consequence of $H_2(M,\z)=\z$.
Indeed, attaching a 1--handle to  $M$ increases
$H_1(M,\z)$ while leaving $H_2(M,\z)$ unchanged.

3. The assumption $H_1(M,\z)=0$ is somewhat arbitrarily chosen.
One could require instead 
$\pi_1(M)=1$.
The variant with $H_1(M,\z)=0$ is the relevant
condition for integral H--cobordism questions
and $\pi_1(M)=1$ connects directly with the 
Montgomery--Yang problem on circle actions.

4. In the algebraic case the inequality 
(\ref{tbmy.thm}.1) holds even if we only assume that
$H_1(M,\q)=0$, but in the topological case this is not enough.
There are lens spaces which bound a rational homology ball,
and taking connected sum of these with $\c\p^2$ gives examples
with an arbitrary number of boundary components.

The simplest algebraic example is the following.
Let $C\subset \c\p^2$ be a smooth conic and $C\subset N$
a regular neighborhood with boundary $L$. 
Since the normal bundle of $C$ in $\c\p^2$ has degree 4,
we see that $L$ is a $\z/4$-quotient of $S^3$.
Set
$M:=\c\p^2\setminus \inter N$. $M$ is a rational homology ball
with $\pi_1(M)=\z/2$ which bounds $L$.

There are many such examples, see, for instance, \cite{ca-ha}.

To get more algebraic ones, let  $u,v$ be  relatively prime
natural numbers.  Then
 the complement of a  regular neighborhood of the curve
$$
C:=(x^{u+v}-yz=0)\subset \p^2(1,u,v).
$$
is a rational homology ball $M$ with $\pi_1(M)=\z/(u+v)$
which bounds a lens space $L$ with $\pi_1(L)=\z/(u+v)^2$.

5. By \cite[p.360]{fi-st-po}, the Poincar\'e sphere
is not H--cobordant to 0, thus one can not use 
connected sums to get counter examples to (\ref{tbmy.thm}).

6. It is possible that Conjecture \ref{tbmy.thm} can be generalized
to arbitrary $H_2(M,\z)$. In this case, 
 (\ref{tbmy.thm}.1),  should be
replaced by 
$$
\sum_i \Bigl( 1-\frac1{|\pi_1(L_i)|}\Bigr)\leq 2+\dim H_2(M,\q).
\eqno{(\ref{tbmy.thm}.2)}
$$
This may hold if $\pi_1(M)=1$ 
 but --  as pointed out to me by P.~Hacking --
$H_1(M,\z)=0$ is not sufficient, not even for algebraic surfaces.
The surfaces in question arise as the minimal compactifications
of Seifert $\c^*$-bundles over $\p^1$.
These were  studied by \cite{dolgachev, pinkham, demazure, fl-za}. 

One of the simplest examples is the following.

\begin{exmp} Rational surfaces $S_m$ with 
$H_1(S_m,\z)=0, H_2(S_m,\z)\cong \z^2$ and with
 $2m$  quotient singularities.
\end{exmp}
  
Let us start with any minimal ruled surface $f:S\to \p^1$.
Pick points $p_i\in F_i\subset S$ where the $F_i$ are different fibers of
$f$ and natural numbers $r_i$. Blow up $(2+r_i)$-times each $p_i\in F_i$
and in each fiber contract all but the last $(-1)$-curve.
We get a singular ruled surface
$$
g=g(r_1,\dots,r_m): S(r_1,\dots,r_m)\to \p^1
\qtq{with} H_2(S(r_1,\dots,r_m),\q)\cong \q^2.
$$
The fiber over $f(p_i)$ is a smooth rational curve with
multiplicity $r_i$. Moreover, $S(r_1,\dots,r_m)$
has 2 singular points  along these fibers.
They are cyclic quotients of the form $\c^2/\frac1{r_i}(1,1)$
and $\c^2/\frac1{r_i}(1,-1)$.

It is easy to see that the fundamental group of the
smooth part is given by generators and relations as:
$$
\pi_1(S(r_1,\dots,r_m)^0)\cong \langle a_1,\dots, a_m: 
a_1^{r_1}=\cdots=a_m^{r_m}=
a_1\cdots a_m=1\rangle.
$$
If the $r_i$ are pairwise relatively prime
then $H_2(S(r_1,\dots,r_m),\z)\cong \z^2$
but the fundamental group is never trivial for $m\geq 3$.
\medskip

 It is worth noting that
(\ref{tbmy.thm}) completely fails 
if $M$ is only a topological manifold.
The number of boundary components can be arbitrary,
see (\ref{top-different}).

\end{say}

\section{The Montgomery--Yang problem}

Fixed point free differentiable circle actions on
$S^7$ with finitely many
nonfree orbits are classified in \cite{mo-ya}, see also \cite{pet}.
Such actions are frequently called {\it pseudo-free}.
The main idea of their classification is the 
following.

We work on $S^7$ and on the
orbifold quotient $X:=S^7/S^1$.
Let $O_i\subset S^7$ be a nonfree orbit
and $x_i\in X$ the corresponding orbifold point.
$O_i$ bounds a disk in $S^7$
whose image is a 2--sphere  $S_i\subset X$
containing $x_i$. 
Since $\dim X>4$, we can arrange these $S_i$ to be disjoint.
The classification now has 2 parts.
\begin{enumerate}
\item[i)] Describe  the
circle action  over a neighborhood of each $S_i$.
These turn out to be diffeomorphic to
linear circle actions on $S^3\times D^4$.
\item[ii)] Describe how these  local models glue together.
\end{enumerate}

By contrast, if we work on $S^5$, then
$\dim S^5/S^1=4$, thus the different $S_i$
do intersect, and we can not
separate the local models from each other. 
(Note that if  a manifold $M$ admits a 
fixed point free differentiable circle action, then
$e(M)=0$, thus this theory is
most interesting for odd dimensional manifolds.)

Montgomery and Yang found that it is
difficult to construct examples
of differentiable circle actions on $S^5$
with many nonfree orbits.
This led to the following:

\begin{conj}[Montgomery--Yang problem] 
\cite[p.41]{mo-ya}
Let $S^1\times S^5\to S^5$ be a differentiable circle action
with only finitely many nonfree  orbits.
Then there are at most 3 nonfree orbits.
\end{conj}

(Note that the proposed partial solution
in \cite{fi-st} is incorrect, see the review, MR0874031.)

There are many different actions
with 3 exceptional orbits.
The simplest ones are linear actions but  there are many
nonlinear examples too.

\begin{exmp} Let $S^5=(|x|^2+|y|^2+|z|^2=1)\subset \c^3$
be the unit sphere. Let $a,b,c\geq 2$ be pairwise relatively prime
natural numbers. Then
$S^1\times S^5\to S^5$, given by 
$$
(\lambda,x,y,z)\mapsto (\lambda^{bc}x,\lambda^{ca}y,\lambda^{ab}z)
 $$
is a differentiable circle action
with  3 nonfree orbits on the 3 coordinate axes.

Note that in this case the quotient
$S^5/S^1$ can be thought of 
as the (complex) weighted projective plane 
$\p^2(a,b,c)$.
\end{exmp}

In general, one can try to classify
circle actions on 5--manifolds $L$ in terms of the
4--dimensional quotient orbifold $L/S^1$.

Let $L$ be a 5--manifold and
$S^1\times L\to L$ a differentiable circle action.
Assume that there are finitely many $S^1$-orbits
$O_1,\dots,O_m$ such that the action is 
free on $L\setminus \cup_iO_i$.
For each $O_i\subset L$ take an $S^1$-invariant
tubular neighborhood 
$O_i\subset T_i\subset L$. Then
$L^0:=L\setminus \cup_iT_i$ is a compact
5--manifold with boundary and with a free circle action.
Thus $M:=L^0/S^1$ is a compact
4--manifold with boundary, where the
boundary components are lens spaces.

The 4--manifold $M$ uniquely determines
$L$ and the $S^1$-action up to diffeomorphism in many cases.
See \cite{mo-ya, fi-st-po} for the pseudo-free case and \cite{e-surf}
in general. (While \cite{e-surf} considers only the algebraic case,
the result is valid in general using the methods
of \cite{circ-act}.) One of the simplest cases is
pseudo-free  circle actions on
5--dimensional rational homology spheres.

\begin{thm}\label{S1.act.char}
There is a one--to--one correspondence between:
\begin{enumerate}
\item Pseudo-free differentiable circle actions on
5--dimensional rational homology spheres $L$  with $H_1(L,\z)=0$.
\item Smooth, compact 4--manifolds $M$  with
boundary such that
\begin{enumerate}
\item  $\partial M=\cup_i L_i$ 
is a disjoint union of lens spaces
$L_i=S^3/\z_{m_i}$,
\item the $m_i$ are relatively prime to each other,
\item $H_1(M,\z)=0$ and $H_2(M,\z)\cong\z$.
\end{enumerate}
\end{enumerate}
Furthermore, $L$ is diffeomorphic to $S^5$ iff
$\pi_1(M)=1$.
\end{thm}

The smooth Bogomolov--Miyaoka--Yau--type conjecture
\ref{tbmy.thm} would give the following 
for circle actions:

\begin{conj} \label{gen.my.conj} Let $L$ be a 
5--dimensional rational  homology sphere with $H_1(L,\z)=0$ admitting
 a pseudo-free differentiable circle action.
Let 
$O_1,\dots,O_k$ be the nonfree orbits with
stabilizers $\z/m_1,\dots,\z/m_k$. Then
$$
\sum_i \bigl(1-\tfrac1{m_i}\bigr)\leq 3.
$$
\end{conj}

In fact, using \cite{circ-act}, one can generalize
this to fixed point free actions which are not
pseudo-free.
Let us say that an orbit $O\subset L$ of an $S^1$-action is {\it exceptional}
if the order of its stabilizer $|\stab(O)|$
is bigger than the least common multiple of the
 orders of the stabilizers of nearby orbits,
that is   $\lcm\{|\stab(O_p)|:p\not\in O\}$.
It is easy to see that the quotient 
$$
m(O):=|\stab(O)|/\lcm\{|\stab(O_p)|:p\not\in O\}
$$
is an integer which is also the order of the local
fundamental group of $L/S^1$ at the image of $O$ \cite[Prop.15]{circ-act}.
Thus (\ref{tbmy.thm}) would imply the following
generalization of (\ref{gen.my.conj}):

\begin{conj} \label{gengen.my.conj} Let $L$ be a 
5--dimensional rational homology sphere with $H_1(L,\z)=0$ admitting
 a fixed point free differentiable circle action.
Let 
$O_1,\dots,O_k$ be the exceptional orbits. Then
$$
\sum_i \bigl(1-\tfrac1{m(O_i)}\bigr)\leq 3.
$$
\end{conj}

\begin{say}
A positive answer to Conjecture \ref{tbmy.thm}
would come  close to settling the
Montgomery--Yang problem. Indeed, 
first we   enumerate
all sequences of pairwise  relatively prime
natural numbers $m_1,\dots,m_k$ such that
$$
\sum_{i=1}^k \bigl(1-\tfrac1{m_i}\bigr)\leq 3.
$$
 The list turns out to be very short and we get 
one of the following cases:
\begin{enumerate}
\item  $k\leq 3$, as required by the Montgomery--Yang problem,
\item $(2,3,5,n)$ for any $(n,30)=1$,
\item $(2,3,7,n)$ for $n\in\{11,13,17,19,23,25,29,31,37,41\}$, or
\item $(2,3,11,13)$.
\end{enumerate}
So, right away we get that there are at most 4
nonfree orbits (assuming Conjecture \ref{tbmy.thm}).

I do not know any 4--manifolds $M$ 
with $H_1(M,\z)=0$ corresponding
to the cases (2--4) above
but I can not exclude these even in the algebraic case.
(See, however, (\ref{4pt.icosa}) for some examples where one
of the singularities is not a cyclic quotient.)

Even if such a manifold exists,
it leads to a counterexample to the
 Montgomery--Yang problem only if it is  simply connected.
\end{say}

The methods to prove (\ref{S1.act.char}) lead to a complete
characterization of  compact, simply connected  5--manifolds
which admit a fixed point free differentiable circle action.
Before we state the result, we need
to define an invariant for   5--manifolds.

\begin{defn}\label{i(L).defn}
For any manifold $L$, write its  second homology as
 a direct sum
of cyclic groups of prime power order
$$
H_2(L,\z)=\z^k+\sum_{p,i} \bigl(\z/p^i\bigr)^{c(p^i)}
\qtq{for some $k=k(L), c(p^i)=c(p^i,L)$.}
\eqno{(\ref{i(L).defn}.1)}
$$
The numbers $k, c(p^i)$ are  uniquely determined by
$H_2(L,\z)$.
One can choose the decomposition
 such that 
the second Stiefel--Whitney class map
$$
w_2:H_2(L,\z)\to \z/2
$$ is zero 
on all but one summand  $\z/2^n$. This value $n$ is unique
and it is   denoted by $i(L)$ \cite{barden}.
Alternatively, 
$i(L)$ is the
smallest $n$ such that there is an $\alpha\in H_2(L,\z)$ such that
$w_2(\alpha)\neq 0$ and $\alpha$ has order $2^n$.

By  \cite{smale, barden}, a compact, simply connected  5--manifolds
is determined by the invariants $k(L), c(p^i,L)$ and $i(L)$.
\end{defn}

\begin{thm}\cite[Thm.3]{circ-act}\label{main.thm}
Let $L$ be a  compact, simply connected  5--manifold
with invariants  $k(L), c(p^i,L), i(L)$ as in (\ref{i(L).defn}.1).
Then $L$ 
 admits a fixed point free differentiable circle action
iff the following conditions hold:
\begin{enumerate}
\item  for every prime $p$, there are
at most $k+1$ nonzero $c(p^i)$,
\item  $i(L)\in \{0,1,\infty\}$, and
\item if $i(L)=\infty$ then there are
at most $k$ nonzero $c(2^i)$.
\end{enumerate}
\end{thm}

\begin{say}[Comments]

1. The above result gives a complete
characterization of those manifolds that admit a 
fixed point free circle action,
but it does {\em not} describe all possible circle actions.
For any $L$ satisfying the above assumptions, infinitely
many topologically distinct circle actions are
constructed in \cite{circ-act} and 
a complete classification seems unlikely.

2. The circle actions constructed in \cite{circ-act}
have 2--parameter families of nonfree orbits.
Pseudo actions exist only if $H_2(L,\z)$
is torsion free (cf.\ \cite[Prop.28]{circ-act}).
\end{say}

\section{H-cobordism of 3--manifolds}

\begin{defn} 
Let $\Sigma$ be an integral homology sphere,
that is $H_*(\Sigma,\z)=H_*(S^3,\z)$. 
$\Sigma$ is {\it H--cobordant} to zero if there is a
smooth 4--dimensional homology cell $W$ 
 with boundary $\partial W=\Sigma$.
Similarly, one can study rational homology spheres
which are rationally H--cobordant to zero.

It is an interesting and difficult question to decide
which integral/rational  homology spheres are 
H--cobordant to zero; see \cite[Chap.7]{sav}
for a  survey.

The case when  $\Sigma$ is Seifert fibered is
very much connected with the 
Bogomolov--Miyaoka--Yau inequality.
\end{defn}

\begin{defn}[Seifert fiber spaces]\label{seifert.def}{\ }

 A {\it Seifert fibered} 3--manifold is a
proper morphism of a 3--manifold to a surface $f:M\to S$
such that every point $s\in S$ has a neighborhood $s\in D_s\subset S$
such that the pair $f^{-1}( D_s)\to  D_s$ is fiber preserving
homeomorphic to one of the  normal forms $f_{c,d}$  defined 
as follows.

Let $S^1(z),D^2(z)\subset \c$ denote the unit circle (resp.\ open unit disk)
where $(z)$ indicates the name of the coordinate.
For a pair of integers $c,d$ satisfying  $0\leq c<d$ and
$(c,d)=1$, define
$$
f_{c,d}:S^1(u)\times D^2(z)\to D^2\qtq{by}
 f_{c,d}(u,z)=u^cz^d.
$$
$f_{c,d}$ restricts to a fiber bundle
$S^1\times (D^2\setminus \{0\})\to D^2\setminus \{0\}$.
The fiber of $f_{c,d}$ over the origin is still $S^1$, but 
$f_{c,d}^{-1}(0)$ has
multiplicity $d$. 

A Seifert fibered 3--manifold $M\to F$ is determined
by the location of the singular fibers $p_i\in F$, 
the above numbers $(c_i,d_i)$ at each point $p_i$
and by a global invariant,
see \cite{seif, seif-book, scott}.
For us the following consequence of the
classification is more important:
\end{defn}

\begin{prop} Given pairwise relatively prime natural numbers
$d_1,\dots,d_k$, there is a unique
integral  homology sphere $\Sigma(d_1,\dots,d_k)$
which has a Seifert fibered structure over $S^2$
with fiber multiplicities $d_1,\dots,d_k$.
\end{prop}

\begin{say}[Seifert disk bundles]
Let $\epsilon$ be a fixed $d$th root of unity
and $\mu_d\subset \c$ the group of $d$th roots of unity.
Consider the quotient
$$
D^2(x)\times D^2(y)/\mu_d(\epsilon, \epsilon^c),
$$
where we use this shorthand to denote the quotient
of $D^2(x)\times D^2(y)$ be the $\mu_d$-action
$(x,y)\mapsto (\epsilon x, \epsilon^cy)$.
The second projection map $y\to y^d$ descends to
$$
F_{c,d}:D^2(x)\times D^2(y)/\mu_d(\epsilon, \epsilon^c) {\to}
 D^2(y)/\mu_d(\epsilon^c)\cong D^2(y^d).
$$
Restricting to $S^1(x)\times D^2(y)$ we get isomorphisms
$$
S^1(x)\times D^2(y)/\mu_d(\epsilon, \epsilon^c)\cong
S^1(x)\times D^2(x^{-c}y)/\mu_d(\epsilon, 1)\cong
S^1(x^d)\times D^2(x^{-c}y),
$$
and setting $u=x^d, z=x^{-c}y$ the projection map $F_{c,d}$ becomes
$$
(u,v)\mapsto u^cz^d=(x^d)^c(x^{-c}y)^d=y^d.
$$

Thus we conclude  that $f_{c,d}$ can be extended 
to the projection map $F_{c,d}$. 

We can glue these together and obtain that
every Seifert fibered 3--manifold $M\to F$ can be obtained
as the boundary of a Seifert disk bundle $D(M)$.
The multiple fibers of multiplicity $m$
correspond to  cyclic quotient singularities
$\c^2/\mu_d(\epsilon, \epsilon^c)$ on $D(M)$.

Note that the zero section gives an embedding $F\into D(M)$
and $D(M)$ retracts to $F$.
The description of $M$ can also be given in terms of 
the self intersection of $F\subset D(M)$.
\end{say}

\begin{lem}\label{seif.h_1.lem}
 Let $M\to F=S^2$ be a Seifert fibered 3--manifold
with pairwise relatively prime  fiber multiplicities $d_1,\dots,d_k$
and corresponding Seifert disk bundle $D(M)$.
Then 
$$
(F\cdot F)=\pm \frac{|H_1(M,\z)|}{d_1\cdots d_k}.
$$
\end{lem}

There are  several series of 
 examples of Seifert fibered  homology spheres
which are H--cobordant to zero
\cite{ak-ki, ca-ha, st-br, neu}, but for all of them
the number of singular fibers is $\leq 3$. 
(See also \cite{orevkov}.)
Let us make the lack of
known examples into a formal conjecture:

\begin{conj}\label{seif.conj}  
Let $M \to S^2$ be a 
 Seifert fibered rational homology sphere
with fiber multiplicities $d_1,\dots,d_k$.
\begin{enumerate}
\item If $M$ is rationally H--cobordant to zero
then $\sum (1-\tfrac1{d_i})<3$.
\item If $M=\Sigma(d_1,\dots,d_k)$ and $M$ is   H--cobordant to zero
then $k\leq 3$.
\end{enumerate}
\end{conj}

It was observed in \cite{fi-st} that this problem
is closely  related to the existence of
4--manifolds whose boundary components are lens spaces,
and hence to Conjecture \ref{tbmy.thm}.

Assume that $\Sigma(d_1,\dots,d_k)$ is the boundary of a
homology cell $W$. Above we wrote
$\Sigma(d_1,\dots,d_k)$ as the boundary of  a
Seifert disk bundle $D(\Sigma(d_1,\dots,d_k))$.
Gluing them together, we get
a 4--dimensional orbifold
$S$. The Mayer--Vietoris sequence shows that
$H_*(S,\z)=H_*(\c\p^2,\z)$.
Thus (\ref{tbmy.thm}) implies that we have at most 4 singular points.

Note, however, that (\ref{seif.conj}.2) is not equivalent
to (\ref{tbmy.thm}). If $S$ is a 4--dimensional orbifold
then usually one can not find an embedded copy of
$S^2$ passing through all singular points 
whose regular neighborhood is  a Seifert disk bundle.

\begin{exmp}[H--cobordisms of lens spaces]\label{lens.H.cob}

We give 2 algebraic constructions  showing that the
lens spaces
$$
\begin{array}{rcl}
L_{(n^2, nc-1)}&:=&S^3/(\z/n^2)(1,nc-1)\\
&:=&(|x|^2+|y|^2=1)/
(x,y)\mapsto (\epsilon x, \epsilon^{nc-1} y)
\qtq{where}\epsilon^{n^2}=1,
\end{array}
$$
are rationally H--cobordant to 0 for any $(c,n)=1$.

The first construction uses deformation of singularities as in
 \cite{wahl2} or \cite[Sec.3]{k-sb}.
Under the action
$$
(x,y)\mapsto (\epsilon x, \epsilon^{nc-1} y)
\qtq{where}\epsilon^{n^2}=1,
$$
the subgroup of $n$th roots of unity acts
as 
$(x,y)\mapsto (\epsilon^n x, \epsilon^{-n} y)$,
 and the quotient is
$$
(uv-w^n=0)\subset \c^3\qtq{where} u=x^n, v=y^n, w=xy.
$$
The induced quotient action of $\z/n\cong (Z/n^2)/(\z/n)$
is given by
$$
(u,v,w)\mapsto (\epsilon^n u, \epsilon^{-n} v, \epsilon^{nc} w)
$$
The Euler characteristic of the Milnor fiber
$(uv-w^n=1)$ is $n$, hence the 
Euler characteristic of the quotient Milnor fiber
$(uv-w^n=1)/(\z/n)$ is $1$. Therefore it is a
rational homology ball.

The second construction uses the  curve
$$
C:=(xy=z^{a+b})\subset \p^2(a,b,1).
$$
$C$ is smooth, rational and passes through 2 quotient
singularities. Thus the bounday $M(a,b)$ of its tubular neighborhood
is Seifert fibered over $C$ with 2 multiple fibers, hence
it is a lens space. Since $(C\cdot C)=(a+b)^2/ab$, we see from
(\ref{seif.h_1.lem}) that $|H_1(L,\z)|=(a+b)^2$.

One needs some explicit computations to decide which lens space.
The chart $z\neq 0$ is isomorphic to $\c^2$
with coordinates $X:=xz^{-a}, Y:=yz^{-b}$. The affine equation
of $C$ is $(XY=1)$.

$\p^2(a,b,1)\setminus (x=0)\cup C$ is isomorphic to
$\c^*\times \c$ and one can choose the isomorphism such that
in the $X,Y$-coordinates the corresponding coordinate functions
are
$$
(XY-1)^aX^{-a-b}\qtq{and} X^{b'-a'}(XY-1)^{-b'}
$$
where $aa'+bb'=1$. Correspondingly, 
$\p^2(a,b,1)\setminus (y=0)\cup C$ is also isomorphic to
$\c^*\times \c$ and  the corresponding coordinate functions
are
$$
(XY-1)^bY^{-a-b}\qtq{and} Y^{a'-b'}(XY-1)^{-a'}.
$$
These determine how the two charts are patched together
along 
$$
\c^*\times \c^*\cong \p^2(a,b,1)\setminus (xy=0)\cup C
$$
with coordinate functions
$$
XY(XY-1)^{-1}\qtq{and} X^bY^{-a}.
$$
A somewhat messy explicit computation gives that 
one gets the lens space
$$
M(a,b)\cong L_{(a+b)^2, (a'-b')(a+b)-1}.
$$
Since $aa'+bb'=1$, we get that
$(a+b)a'-(a'-b')b=1$, thus $a+b$ and $a'-b'$ are relatively prime.

It is easier to see this by putting together the
 two descriptions  as follows.
Fix $n$ and $a,c$ such that $ac\equiv 1\mod n$.
Consider the family of weighted affine hypersurfaces
$$
X(\lambda):=(xy-z^n+\lambda t=0)\subset \p(a,n-a,1,n)\setminus (t=0).
$$
For $\lambda=0$ we get 
$$
X(0)=(xy-z^n=0)/\tfrac1{n}(a,n-a,1)
\cong (xy-z^n=0)/\tfrac1{n}(1,-1,c).
$$
For $\lambda\neq 0$ we can eliminate $t$ to get that
$$
X(\lambda)\cong \p(a,n-a,1)\setminus (xy-z^n=0).
$$
Thus the quotient Milnor fiber of the first constructions
is isomorphic to the complement of the curve $C$
in the second construction.

\end{exmp}

The following series of 
\cite[5.9.2]{wahl-sm}
gives an algebraic realization of
the rational H--cobordisms constructed in \cite[p.132]{neu}.
See  (\ref{qs.ratcurves}) for its projective version.

\begin{exmp}\label{H.cob.exmps}
 Let $a,b,c\geq 1$ be  integers. Set
$S_{abc}:=(x^ay+y^bz+z^cx=0)\subset\c^3$.
$S_{abc}$  is quasi homogeneous with weights
$(bc-c+1, ca-a+1, ab-b+1)$,
thus its link is a Seifert fibered 3--manifold $L_{abc}\to S^2$
with 3 multiple fibers of multiplicities
$ab-b+1, bc-c+1, ca-a+1$.
The Milnor fiber is
$$
M_{abc}:=(x^ay+y^bz+z^cx=1)\subset\c^3,
$$
and its Euler characteristic  is
$abc+1$. Let $\omega$ be a primitive $(abc+1)$st root of unity.
The cyclic group $\z/(abc+1)$ acts freely on $M_{abc}$ by
$$
(x,y,z)\mapsto (\omega x, \omega^{-a}y,\omega^{ab}z).
$$
The quotient $M_{abc}/\bigl(\z/(abc+1)\bigr)$
has Euler characteristic 1, thus it is a
rational homology ball. 
It bounds the Seifert fibered 3--manifold $L_{abc}/\bigl(\z/(abc+1)\bigr)$
which also has  3 multiple fibers of multiplicities
$ab-b+1, bc-c+1, ca-a+1$.
\end{exmp}

The following remarkable example of  \cite{w-let} shows that there are 
even examples with 4 singular fibers. Further such examples are in 
\cite{st-sz-w}.

\begin{exmp}[J.\ Wahl] A Seifert fibered rational homology sphere with
4 multiple fibers which is rationally H--cobordant to 0.

The classification of all finite subgroups of $O(n)$
which act freely on the sphere $S^{n-1}$ is the
 spherical space form problem.
 See  \cite[Part III]{wolf} for a thorough treatment.
One such action is obtained as follows
(cf.\ \cite[5.5.6]{wolf}).

Fix an integer $u\geq 2$ and set
$$
m:= 3u^2-3u+1,\ n:=9u,\ r:=3u^2-6u+2,\ a:=3u-1.
$$
Let $\zeta$ (resp.\ $\omega$)  be a primitive $m$th 
(resp.\ $3u$th) root of unity
and consider the subgroup $G\subset GL(3,\c)$ 
generated by
$$
A:= (x,y,z)\mapsto (\zeta x,\zeta^r y,\zeta^{r^2}z)
\qtq{and}
B:= (x,y,z)\mapsto (\omega^{-1}z, x,y).
$$
$G$ is nilpotent, 
its order is $nm$ and it is given by
the relations
$$
\langle A,B: A^m=B^n=1, BAB^{-1}=A^r\rangle.
$$
Consider now the function
$$
f(x,y,z):= x^ay+y^az+\omega z^ax.
$$
\medskip
{\it Claim.} Notation as above.
\begin{enumerate}
\item  The link of the singularity
$(f=0)/G$ is  Seifert fibered with 4 singular fibers
of multiplicities  $(3,3,3,m)$.
\item The Milnor fiber $(f=1)/G$ is a rational homology ball.
\end{enumerate}
\medskip

Proof. The singularity $(f=0)$ 
can be resolved by 1 blow up. The exceptional curve
$C$ is  smooth  of genus $\binom{a}{2}$.

Note that $B^3=\omega{\mathbf 1}$ acts trivially on $C$.
The group $G/\langle B^3\rangle$ acts on $C$ with 4
nonfree orbits. These are:

-- the orbits of $(1:\lambda \omega : \lambda^2 \omega)$
where $\lambda^3 \omega=1$; these have order 3 stabilizers, and

-- the orbit of $(1:0:0)$; this has order $m$ stabilizer.

\noindent The Hurwitz formula now shows that
$C/G\cong \p^1$. Moreover $B_0\c^3/G$ has 4 singular points,
three have index 3 and one has index $m$. This shows the first part.

To see the second part, note that
by explicit computation,
the  Milnor number of $f$ is 
$$
\dim \c[x,y,z]/(\partial f/\partial x,
\partial f/\partial y, \partial f/\partial z)= a^3.
$$
Thus the Euler characteristic of the Milnor fiber $(f=1)$ is
$1+a^3=mn=|G|$. Therefore
 the Euler characteristic of the quotient Milnor fiber
  $(f=1)/G$ 
is 1, and hence it is a rational homology ball.\qed
\end{exmp}

\begin{say}[Smoothing surface singularities]
Let $(0\in S)$ be a normal surface singularity
with link $L$ and let $M$ be the Milnor fiber of a smoothing of
$S$. With suitable care, $M$ is a 4--manifold whose boundary is
$L$. We are especially interested in the case when
$M$ is a rational homology ball.

This imposes a very strong restriction on $(0\in S)$.
For instance,  it seems to have been known  for some time
 that $(0\in S)$ is a rational singularity,
cf.\ \cite[Sec.2.3]{st-sz-w}. 
The latter also shows that (\ref{seif.conj}.1) holds
for the link of a normal surface singularity
 which has a smoothing 
whose Milnor fiber is a rational homology ball.

For other relationships between the algebraic geometry
of a surface singularity and the Seiberg-Witten theory
of its link see \cite{ne-wa, w-char, n-n1, n-n2, n-n3}.
\end{say}

\begin{rem}\label{top-different} By the results of 
\cite{free}, every  integral homology sphere
bounds a topological homology 4--cell.
By the above construction this implies that
there are topological 
integral 
homology $\c\p^2$-s with any number of singular points.
\end{rem}

\section{Open problems on algebraic surfaces}

On the algebraic geometry side, all of the
questions can be gathered into one central problem:

\begin{prob}\label{hopeless} Classify all integral/rational
 homology $\c\p^2$-s with
quotient singularities.
\end{prob}

I do not expect this to be  feasible.
The case when the canonical class is
anti ample, that is, we are looking at log--Del Pezzo
surfaces, received a lot of attention,  but
there seem to be too many cases for a complete structure theorem;
see \cite{miyanishi, ke-mc, shokurov}.

The case where the canonical class is numerically trivial should be
hard but there is a clear path to follow.
If $mK_S\sim 0$ then $S$ has a degree $m$ cover
which is either an Abelian surface or a K3 surface with
Du~Val singularities.  Thus the problem reduces to the
classification of cyclic group actions on
Abelian and K3 surfaces.
The 7 cases where the Picard number of
the K3 surface is maximal  are classified in \cite{og-zh}.

Very little is known about the case when the canonical class is
ample. The recent classification in the smooth case
\cite{prasad} is very significant, but it probably
says very little about the singular case.

The following examples are worked out in (\ref{rhcp2.say}).

\begin{exmp} \label{rat.surf.exmp} Let $a_1,a_2,a_3,a_4$ be
natural numbers  such that
$a_2a_3a_4-a_3a_4+a_4-1$ and  
$a_1a_2a_3a_4-1$ are relatively prime.
(This holds in at least 75\% of all cases.)
Let $S$ be the surface
$$
S:=S(a_1,\dots,a_4):=
(x_1^{a_1}x_2+x_2^{a_2}x_3 + x_3^{a_3}x_4+ x_4^{a_4}x_1=0)
\subset \p(w_1,\dots,w_4)
$$
where, using subscripts modulo 4,
$$
w_i=a_{i+1}a_{i+2}a_{i+3}-a_{i+2}a_{i+3}+a_{i+3}-1.
$$
Note that  $S$ contains the curves
$(x_1=x_3=0)$ and $(x_2=x_4=0)$ and they can be 
contracted. Thus we obtain  a surface
 $S^*=S^*(a_1,a_2,a_3,a_4)$. 
We compute in Section 5 that:
\begin{enumerate}
\item $S^*$ is a rational surface with
quotient singularities,
\item  $H_*(S^*,\q)=H_*(\c\p^2,\q)$,
\item if $a_i\geq 4$ then 
the canonical class of $S^*$ is ample, and
\item $(K_{S^*}^2)$ converges to 1 as $\min \{a_i\}\to \infty$.
\end{enumerate}
\end{exmp}

Let us also note that in positive characteristic
one can get examples with many quotient singularities.

\begin{exmp}[Characteristic $p$]  Let $k$ be a  field 
of characteristic $p$. For some $q=p^m$, blow up all the
$\f_q$-points of $\p^2$. The birational transforms of the
$\f_q$-lines become disjoint, smooth, rational curves with self intersection
$-q$. They can be contracted to obtain a surface
$X_q$ defined over $\f_q$ with $q^2+q+1$ singular points.
These are quotient singularities of the type $\a^2/\mu_q$
where $\mu_q$ is the subgroup scheme of $q$th roots of unity
$\spec k[t,t^{-1}]/(t^q-1)$. 
(Note that in characteristic $p$ it is the 
$\z/q$ quotients that behave very badly (cf.\ \cite{artin})
and the $\mu_q$-quotients are the correct characteristic $p$
analogs of characteristic $0$ quotient singularities.)

$X_2$ is a Del Pezzo surface of degree 2, but 
the canonical class is ample for $q\geq 3$. 
\end{exmp}

Below I list some special cases of Problem \ref{hopeless}
which are of interest either as partial steps in the
classification or as having immediate topological
consequences.

\begin{conj}\cite[4.17]{e-surf}\label{sconnconj}
 Let $S$ be a rational homology $\c\p^2$ with
quotient singularities.
If $S^0:=S\setminus \sing S$ is simply connected then $S$ is rational.
\end{conj}

This is similar in spirit  to the result 
that in dimension 2, all algebraic  $\q$--homology cells are rational
\cite{gu-pr}. (\ref{sconnconj}) was verified by \cite{keum}
when the singularities are not very complicated.

\begin{conj}[Algebraic Montgomery--Yang problem]
 Let $S$ be a rational homology $\c\p^2$ with
quotient singularities.
If $S^0:=S\setminus \sing S$
 is simply connected then $S$ has at most 3 singular points.
\end{conj}

\begin{exmp}\label{4pt.icosa}
 Let $G\subset SL(2,\c)$ be 
subgroup which contains no quasi--reflections such that its
image in $PSL(2,\c)$ is the  icosahedral group $I$
(see \cite{briesk-s} for a complete list and the
corresponding quotient singularities). Let $Z\subset G$ be its center,
then $G/Z\cong I$.
Extend the $G$--action on $\c^2$ to $\c\p^2$. 
The center acts trivially on the line at infinity
and $\c\p^2/Z$ is a cone over the rational normal  curve of degree $|Z|$.
Then
$S_G:=\c\p^2/G=(\c\p^2/Z)/I$ has 4 quotient  singularities,
one of type $\c^2/G$ at the origin, 3 of types
$\c^2/\z_2, \c^2/\z_3,\c^2/\z_5$ at infinity. The fundamental group
of $S^0_G$ is $I$, thus $S_G$ is an integral homology $\c\p^2$.

These examples can also be obtained by starting with a minimal ruled surface
and blowing up inside 3 of the fibers.
\end{exmp}

\begin{prob} Classify all integral/rational homology $\c\p^2$-s with
4 or more quotient singularities.
\end{prob}

 It is also of considerable interest to study algebraic surfaces
that lead to H--cobordisms of Seifert fibered manifolds.
If we want to stay completely algebraic, then 
we are lead to the following

\begin{prob}\label{alg.seif.Hcob} Classify all  pairs $(S,C)$ such that
\begin{enumerate}
\item $S$ is a rational homology $\c\p^2$,
\item $C$ is a rational curve, homeomorphic to $S^2$.
\end{enumerate}
\end{prob}

In this case $S\setminus C$ is a {\it rational homology plane}
or a {\it $\q$-acyclic surface}. That is a (nonproper) surface 
$X$  such that $H_*(X,\q)=H_*(\c^2,\q)$.

There is a huge body of literature devoted to classifying
integral/rational homology planes. 
See \cite{1, 5, 2, 3, 4, 6, gu-pr, dr1, dr2, dr3}
 and the many references there.
Nonetheless, most rational homology planes can not be compactified
to get a rational homology $\c\p^2$, thus 
(\ref{alg.seif.Hcob}) may be a much easier problem.

\begin{prob} Topological smoothings of surfaces with
quotient singularities.
\end{prob}\label{top.smooth.prob}

Let $S$ be a proper surface with quotient singularities
$p_i\in S$. Let $p_i\in M_i$ be small conical neighborhoods.
 For each singularity $(p_i\in M_i)$ choose a 
4--manifold such that $\partial N_i=\partial M_i$.
For instance, $N_i$ could be the general fiber of a smoothing
of $M_i$.
 In general, there need not be an
algebraic deformation of $S$ which realizes these local smoothings.
However, one can always get a differentiable manifold
by replacing each singular $M_i$ with the smooth $N_i$.

It is interesting to find this way 
differentiable manifolds which are not algebraic.
This idea, though formulated purely topologically, 
was used in \cite{st-sza, p-s-sz, fi-st-5} to construct 
exotic differentiable structures on $\c\p^2$
blown up at $\geq 5$ points.

It is especially interesting to work with the case
when each $N_i$ is a rational homology ball.
In the early seventies, Casson, Gordon and Conway
(unpublished) 
found 3  such classes:
\begin{enumerate}
\item $\c^2/\frac1{n^2}(1,na- 1)$ where $(n,a)=1$,
\item $\c^2/\frac1{n^2}(1,d(n- 1))$ where $d|n- 1$ is odd, and
\item $\c^2/\frac1{n^2}(1,d(n- 1))$ where $d|2n+ 1$.
\end{enumerate}
Recently \cite{lisca} proved that these are in fact all,
if we also take into account two elementary observations:
\begin{enumerate}\setcounter{enumi}{3}
\item replacing the generator of $\z/p$ by its inverse
shows that $\c^2/\frac1{p}(1,q)\cong \c^2/\frac1{p}(1,q')$
for $qq'\equiv 1 \mod p$, and
\item  conjugating one of the coordinates shows
that $\c^2/\frac1{p}(1,q)$ is diffeomorphic to $\c^2/\frac1{p}(1,p-q)$.
\end{enumerate}

Two algebraic realizations of the first series were given
in (\ref{lens.H.cob}). I do not know algebraic descriptions
of the other two.

Let $S$ be a proper rational homology $\c\p^2$
 with quotient singularities
$p_i\in S$ which are on the list
(\ref{top.smooth.prob}.1--5).
The resulting topological smoothing is then a smooth
4--manifold which is also a  rational homology $\c\p^2$.

In an earlier version of this note I
raised the possibility that such surfaces could lead to
a fake $\c\p^2$, that is, a smooth 4--manifold
homeomorphic but not diffeomorphic to $\c\p^2$.
Moreover,  I was hoping to do this with smoothings that
are locally complex analytic.
To this end, one needs to 
find rational homology $\c\p^2$-s with
singularities of the form  $\c^2/\frac1{n^2}(1,na-1)$.

If $S$ has such  singularities and $-K_S$ is ample,
then the local deformations of the singular points can be globalized,
and $S$ is the degeneration of smooth Del Pezzo surfaces.
These were studied and partially classified in
\cite{manetti, hac-pro}.
The  general type examples, however,  have  a tendency to be
rigid (cf.\ \cite[6.12]{6}), and they may lead to
new differentiable 4--manifolds.
Unfortunately, we get nothing interesting.

\begin{thm}[P.\ Hacking] Let $S$ be a 
rational homology $\c\p^2$ with
singularities of the form  $\c^2/\frac1{n_i^2}(1,n_ia_i-1)$
for some $(n_i,a_i)=1$.
If $K_S$ is nef then $S$ is smooth.
\end{thm}

Proof. Since $K_S$  is nef,
we have the orbifold BMY inequality 
$c_1(S)^2 \le 3e_{orb}(S)$. 

Noether's formula $\chi(\o_S)=\frac1{12}(c_1(S)^2+e(S))$
usually needs a correction term for
each quotient singularity, but for our singularities
the
correction term vanishes, cf.\ \cite[Prop.3.5]{hac-pro}.
Thus  $c_1(S)^2=9$ and so $3\leq e_{orb}(S)$ thus $S$ is smooth.\qed
\medskip

If we also allow the other singularities on the
list (\ref{top.smooth.prob}.1--5), there are
interesting examples. For instance
$$
\bigl(u^2=(x^2+y^2-z^2)(x^2+2y^2-z^2)\bigr)\subset \p^3(1,1,1,2)
$$
is a degree 2 Del Pezzo surface with Picard number 2
and two singular points of type $A_3$, that is,
$\c^2/\frac14(1,3)$. Topological smoothing
creates out of it a 4--manifold with $b_2=2$
(which is probably not simply connected).

It would be interesting to start a systematic study of such examples.

\section{Examples of rational homology projective spaces}

In this section we investigate hypersurfaces in weighted projective
spaces given by an equation
$$
H(a_1,\dots,a_n):=
(x_1^{a_1}x_2+x_2^{a_2}x_3+ \cdots + x_{n-1}^{a_{n-1}}x_n+ x_n^{a_n}x_1=0)
\subset \p(w_1,\dots,w_n).
$$
These  hypersurfaces,
or rather, the corresponding cones, play a fundamental role
in the classification of quasi--homogeneous singularities
 \cite{or-ra, kou}, but they  have
many other remarkable properties as well.

The best known examples arise when
$a_1=\cdots=a_n=a$, giving hypersurfaces in
ordinary projective space
$$
H(a):=
(x_1^{a}x_2+x_2^{a}x_3+ \cdots + x_{n-1}^{a}x_n+ x_n^{a}x_1=0)
\subset \p^{n-1}.
$$
These have been studied for their  large group of automorphisms
among others. (The case $a=n=3$ is Klein's curve of genus 3
with  a simple group of order 168 as automorphisms.)

It turns out, however, that this case is  completely misleading
 and for general $a_1,\dots,a_n$ we get
very different behavior.

The closely related examples of 
\cite{libg, bar-dim}
$$
(x_1^{a}x_2^{d-a}+x_2x_3^{d-1}+ \cdots + x_{n-1}x_n^{d-1}+ x_{n+1}^{d}=0)
\subset \p^n
$$
give hypersurfaces with (non--quotient) isolated singularities
which have the same integral homology  as $\c\p^{n-1}$
if $\gcd(a,d(d-1))=1$.

\begin{say}[Summary of the results]\label{summary.say}

Let $H=H(a_1,\dots,a_n)$ be as above.
Then $H$
has only cyclic quotient singularities for $n\geq 4$.
Under a mild but not very explicit restriction on $ a_1,\dots,a_n$
(\ref{num.conds}) we show that
\begin{enumerate}
\item $H$ is birational to $\c\p^{n-2}$, but
\item if every $a_i\geq n$ and $n\geq 4$
then the canonical class $K_{H}$ is ample and
its self intersection $(K_H^{n-2})$ converges to 1 
 as $\min\{ a_i\}\to \infty$.
\end{enumerate}
Moreover, if $n$ is odd then 
\begin{enumerate}\setcounter{enumi}{2}
\item  $H$ is a  rational homology  $\c\p^{n-2}$, and
\item $\p(w_1,\dots,w_n)\setminus H$ is a rational homology $\c^{n-1}$.
\end{enumerate}
For $n=3$ this gives many examples of 
Seifert fibered rational homology spheres which are
rationally $H$-cobordant to 0, see (\ref{H.cob.exmps}).

If $n=2m$ is even then  $H_{n-2}(H,\q)$ has dimension 3. 
However, if $n=4$,  then $H$ contains
2 disjoint contractible curves and after
contracting them we get $H\to H^*$ and in (\ref{rat.surf.exmp})
we show that
\begin{enumerate}\setcounter{enumi}{5}
\item  $H^*$ is a  rational homology  $\c\p^{n-2}$, and 
\item if every $a_i\geq 4$ then  the canonical class $K_{H^*}$ is ample.
\end{enumerate}
\end{say}

\begin{rem} In singularity theory and in topology the
Brieskorn--Pham singularities
$$
(x_1^{a_1}+x_2^{a_2}+ \cdots + x_n^{a_n}=0)\subset \c^n
$$
are much better known.
When the link of a Brieskorn--Pham singularity
is a homology sphere, then the corresponding projective hypersurface
is isomorphic to a weighted projective space \cite{briesk}
and all the intricate geometry is concentrated
in the corresponding Seifert bundle structure 
(see, for instance, \cite{or-wa}).

By contrast, for the singularities
$$
(x_1^{a_1}x_2+x_2^{a_2}x_3+ \cdots + x_{n-1}^{a_{n-1}}x_n+ x_n^{a_n}x_1=0)
\subset \c^n
$$
the projective hypersurfaces are also very interesting
and the Seifert bundle structure  is usually simple.
\end{rem}

\begin{say}[Numerical conditions]\label{num.conds}

The exponents $a_1,\dots,a_n$,
the weights $w_1,\dots,w_n$ and the degree $d$ are related by the
equations
$$
a_iw_i+w_{i+1}=d\qtq{for} i=1,\dots,n,
$$
 where we write all subscripts modulo $n$.
Let us fix the exponents $a_i$ and assume that $\prod a_i\neq (-1)^n$.
Using (\ref{linalg.lem}) the system can be solved
explicitly. Set
$$
W_i:=\sum_{j=1}^n(-1)^{j-1}\prod_{\ell=i+j}^{i+n-1}a_{\ell}
\qtq{and} D:=\prod_{\ell=1}^n a_{\ell}+(-1)^{n-1},
$$
where $D$ is the determinant of the system.
It is easy to check that
$$
a_iW_i+W_{i+1}=D\qtq{for} i=1,\dots,n.\eqno{(\ref{num.conds}.1)}
$$
We would like to have a well formed weighted projective space,
thus we have  to divide the weights by
their greatest common divisor
$$
w^*:=\gcd(W_1,\dots,W_n).\eqno{(\ref{num.conds}.2)}
$$
Note that the equations (\ref{num.conds}.1) imply that
$$
w^*=\gcd(W_i,D)\ \forall i\qtq{and}
w^*=\gcd(W_i,W_{i+1})\ \forall i.\eqno{(\ref{num.conds}.3)}
$$
It turns out that the cases with $w^*=1$
have many special properties that are not shared by the
examples with $w^*>1$. 

 It is not clear to me how to determine whether $w^*=1$,
other than actually computing it. It is, however, easy to see
that $w^*=1$ happens frequently.

Note that we can write $W_n=a_{n-1}A\pm 1$ and
$D=a_na_{n-1}B\pm 1$ where $A,B$ depend on
$a_1,\dots,a_{n-2}$ only. Fix a prime $p$
and $a_1,\dots,a_{n-2}$.
No matter what $A,B$ are, there is at most one choice
for $a_{n-1}$ modulo $p$ such that 
 $p$ divides $W_n$ and then at most one choice
for $a_{n}$ modulo $p$ such that 
 $p$ divides $D$.
 By the Chinese remainder theorem,
the conditions for different primes are independent.
Thus the proportion of the $n$-tuples  with $w^*=1$ is
asymptotically at least
$$
\prod_p\bigl( 1-\frac{1}{p^2}\bigr)=\frac{1}{\zeta(2)}=\frac {6}{\pi^2}
=0.607...
$$
Note that $D=\prod a_i\pm 1$ is divisible by 2 only if all the
$a_i$ are odd, and then $w_i\equiv n\mod 2$. This 
allows us to conclude that the above $0.607$ can be
improved to $0.8$ if $n$ is odd and to $0.75$ if $n$ is even.

Set
$$
w_i:=\tfrac1{w^*}W_i\qtq{and} 
 d:= \tfrac1{w^*}D.\eqno{(\ref{num.conds}.4)}
$$
From (\ref{num.conds}.3) we conclude that
$$
\gcd(w_i,d)=1\ \forall i\qtq{and}
\gcd(w_i,w_{i+1})=1\ \forall i.\eqno{(\ref{num.conds}.5)}
$$
\end{say}

Note  that
$H$  can also be viewed as a general element of the
linear system
$$
|H|=|x_1^{a_1}x_2,x_2^{a_2}x_3, \cdots, x_{n-1}^{a_{n-1}}x_n, x_n^{a_n}x_1|.
$$
Indeed, take any $\lambda_i\neq 0$ and consider 
$$
H(\lambda):=(\lambda_1x_1^{a_1}x_2+\lambda_2x_2^{a_2}x_3+ \cdots + 
\lambda_{n-1}x_{n-1}^{a_{n-1}}x_n+ \lambda_nx_n^{a_n}x_1=0)
$$
Choose $\mu_i$ such that $\mu_i^{a_i}\mu_{i+1}=\lambda_i^{-1}$.
Such a choice is possible since we can view these equations as a 
linear system for the $\log \mu_i$ whose determinant is
$D\neq 0$. Thus $H$ and $H(\lambda)$ differ only by a 
coordinate change.

\begin{thm}\label{main.exmp.thm}
 Assume that $\prod a_i\neq (-1)^n$
and write all subscripts modulo $n$. 
Define $w_i,d$ and $w^*$ as in (\ref{num.conds}.4). Then
\begin{enumerate}
\item $\p(w_0,\dots,w_n)$ is a 
well formed weighted projective space
whose singular set has dimension $\leq \rdown{n/2}-1$.
\item The hypersurface
$$
H(a_1,\dots,a_n):=
(x_1^{a_1}x_2+x_2^{a_2}x_3+ \cdots + x_{n-1}^{a_{n-1}}x_n+ x_n^{a_n}x_1=0)
\subset \p(w_1,\dots,w_n)
$$
is quasi--smooth.
\item $\p(w_1,\dots,w_n)\setminus H(a_1,\dots,a_n)$
is smooth.
\item If $w^*=1$ then $H$ is birational to $\c\p^{n-2}$.
\item $K_H$ is ample if $\min\{a_i\}\geq n$ and 
 $(K_H^{n-2})\to 1$  as $\min\{a_i\}\to \infty$.
\item If $n$ is odd and $w^*=1$ then
$$
\begin{array}{l}
H_*(H(a_1,\dots,a_n),\q)\cong H_*(\c\p^{n-2},\q)\qtq{and}\\
H_*(\p(w_1,\dots,w_n)\setminus H(a_1,\dots,a_n),\q)\cong H_*(\c^{n-1},\q).
\end{array}
$$
\item If $n=2m$ is even and $w^*=1$ then
$$
\begin{array}{l}
H_j(H(a_1,\dots,a_n),\q)\cong H_j(\c\p^{n-2},\q)\qtq{for $j\neq n-2$, and}\\
H_{n-2}(H(a_1,\dots,a_n),\q)\cong \q^3.
\end{array}
$$
The middle homology is spanned by
the the complete intersection class $c_1(\o_{\p}(1))^{m-1}\cap [H]$ and the 
two disjoint weighted linear subspaces
$$
(x_1=x_3=\cdots=x_{2m-1}=0)\qtq{and}
(x_2=x_4=\cdots=x_{2m}=0).
$$
\end{enumerate}
\end{thm}

Proof. The singular locus of
$\p(w_1,\dots,w_n)$ is a union of weighted linear subspaces
$L_I$ where
 $I\subset \{1,\dots,n\}$ is a  subset such that
$\gcd\{w_i:i\in I\}\neq 1$ and
 $$
L_I:=\{(x_1,\dots,x_n):x_j=0\ \forall \ j\not\in I\}.
$$
As we noted in (\ref{num.conds}.5), $I$ does not contain 
any pair of indices whose
difference is 1. Thus $|I|\leq n/2$ and so
$\dim L_I\leq \rdown{n/2}-1$.
We also see that $L_I\subset H$ for every such $I$.
This shows (\ref{main.exmp.thm}.1) and (\ref{main.exmp.thm}.3). 

Outside $(\prod x_i=0)$ the hypersurface
$H$ is smooth by Bertini.
Assume that $H$ is not quasi--smooth at the point
$(p_1,\dots,p_n)$ and $p_i=0$. 
Then $\partial h/\partial x_i=0$ shows that $p_{i-1}=0$
and by repeating the argument we get that all the $p_j=0$. Thus
$H$ is quasi--smooth, proving (\ref{main.exmp.thm}.2).

Assuming that $w^*=1$,  we show that the
linear system
$$
|H|=|x_1^{a_1}x_2,x_2^{a_2}x_3, \cdots, x_{n-1}^{a_{n-1}}x_n, x_n^{a_n}x_1|
$$
maps $\p(w_1,\dots,w_n)$ birationally to $\p^{n-1}$ and so
 $H$ is mapped birationally to
a hyperplane in $\p^{n-1}$. Note that $|H|$ 
restricts to a homomorphism between the tori
$$
\eta:(\c^*)_{\mathbf x}^n\to (\c^*)_{\mathbf y}^n
\qtq{given by} \eta^*y_i=x_i^{a_i}x_{i+1},
$$
where $(y_1,\dots,y_n)$ are coordinates on the target $\p^{n-1}$.
The degree of $\eta$  is the determinant of the
matrix of exponents, which we already computed to be
$D=\prod_i a_i+(-1)^{n-1}$.

Let us now restrict $\eta$ to  the 1--parameter subgroup
$(\lambda^{w_1},\dots, \lambda^{w_n})$. We get
a homomorphism of degree $d$:
$$
\eta:(\lambda^{w_1},\dots, \lambda^{w_n})\to 
(\lambda^d,\dots,\lambda^d).
$$
Note that  $w^*=1$ iff $d=D$,
thus if $w^*=1$ then $\eta$ descends to an isomorphism
$$
(\c^*)_{\mathbf x}^n/(\lambda^{w_1},\dots, \lambda^{w_n})\to
 (\c^*)_{\mathbf y}^n/(\lambda,\dots,\lambda).
$$
This is exactly the map given by $|H|$,
proving (\ref{main.exmp.thm}.4).

Although it is not needed for our purposes,
one can also write down the inverse of the map
given by $|H|$. First we get the formulas
$$
\begin{array}{lcl}
x_i^D&=&\prod_{j=0}^{n-1}y_{i+j}^{b_{ij}}\qtq{where}
b_{ij}:=(-1)^{j+1}\prod_{\ell=i+j-1}^{i+n-1}a_{\ell}\qtq{and}\\
x_i^{w_{i+1}}x_{i+1}^{-w_i}&=& x_i^Dy_i^{-w_i}.
\end{array}
$$
Then one can easily check that the monomials
$x_i^D$ and $x_i^{w_{i+1}}x_{i+1}^{-w_i}$
generate the subring $\c(x_1,\dots,x_n)^{(D)}$ of those
elements whose degree is divisible by $D$.
Thus we get an explicit isomorphism
$$
|H|^*:\c(y_1,\dots,y_n)\cong \c(x_1,\dots,x_n)^{(D)}.
$$
Nevertheless, I found it very difficult to compute
anything based on these formulas.

For $n\geq 5$ the singular set of $\p(w_1,\dots,w_n)$
has codimension $\geq 3$, thus the 
 canonical class of $H$ is given by the adjunction formula,
$$
K_H=(K_{\p}+H)|_H=\o_{\p}(d-\sum w_i)|_H.
$$
The minor modifications needed in the
 few cases when $n=4$ and $\p(w_1,\dots,w_n)$
has 1--dimensional singular set are discussed in
(\ref{rhcp2.say}).
If $w^*=1$ then the two highest terms in the coefficient of $K_H$ are
$$
\bigl(\ \prod a_i\bigr) \bigl(1-\sum \tfrac1{a_i}\bigr),
$$
which is positive as soon as $\min\{a_i\}\geq n$. In fact it is easy to see
that $d-\sum w_i>0$ if $\min\{a_i\}\geq n$.
The self intersection of $K_H$ is computed asymptotically by 
$$
(K_H^{n-2})=(\o_{\p}(d-\sum w_i)^{n-2}\cdot \o_{\p}(d))
\sim \frac{d^{n-1}}{\prod w_i}\sim w^*.
$$
 
We compute the homology groups of $H$ using the
Milnor--Orlik formula  \cite{mi-or}.

Let $f(x_1,\dots,x_n)$ be  a weighted homogeneous polynomial
of weighted degree $d$ where the variable $x_i$ has weight $w_i$.
Assume that $(f=0)$ has an isolated singularity at the origin
and let $L=L(f):=(f=0)\cap S^{2n-1}(1)$ be its link.
Then $L$ is $(n-3)$--connected and the rank of the
middle homology groups is given by 
$$
\dim H_{n-2}(L,\q)=
\sum_{I\subset \{1,\dots,n\}} (-1)^{n-|I|}
\frac{\prod_{i\in I} (d/w_i)}{\lcm\{u_i:i\in I\}},
$$
where we write $\frac{d}{w_i}=\frac{u_i}{v_i}$ in lowest terms.

In our case $(d,w_i)=1$ for every $i$  (\ref{num.conds}.5)
an so $u_i=d$. Thus 
$\lcm\{u_i:i\in I\}=d$ save for $I=\emptyset$. Thus
the formula becomes
$$
\begin{array}{l}
\sum_{I\subset \{1,\dots,n\}} (-1)^{n-|I|}
\frac{\prod_{i\in I} (d/w_i)}{\lcm\{u_i:i\in I\}}=\\
\qquad\qquad\qquad\qquad=(-1)^{n}+\frac{(-1)^{n-1}}{d}+
\frac1{d}\sum_{I\subset \{1,\dots,n\}} (-1)^{n-|I|}
\prod_{i\in I} \frac{d}{w_i}\\
\qquad\qquad\qquad\qquad
=(-1)^{n}+\frac{(-1)^{n-1}}{d}+\frac1{d}\prod_i\Bigl(\frac{d}{w_i}-1\Bigr)\\
\qquad\qquad\qquad\qquad
=(-1)^{n}+\frac{(-1)^{n-1}}{d}+
\frac1{d}\prod_i\Bigl(\frac{a_{i-1}w_{i-1}}{w_i}\bigr)\\
\qquad\qquad\qquad\qquad=(-1)^{n}+\frac{(-1)^{n-1}}{d}+\frac1{d}\prod_i a_i\\
\qquad\qquad\qquad\qquad=(-1)^n+w^*,
\end{array}
$$
where at the last step we took into account that 
$dw^*=D=\prod_i a_i +(-1)^{n-1}$.

The link $L$ is a Seifert $S^1$-bundle
over $X:=(f=0)\subset \p(w_1,\dots,w_n)$
and the resulting Leray spectral sequence
is  easy to compute (with rational coefficients),
see, e.g. \cite{or-wa}. This gives
(\ref{main.exmp.thm}.6--7) except for the
precise identification of
$H_{n-2}(H,\q)$ in the $n=2m$ case.

We use this only for $n=4$, where it is worked out in
(\ref{rhcp2.say}).\qed

\begin{lem}\label{linalg.lem} Assume that $\prod a_i\neq (-1)^n$
and write all subscripts modulo $n$. 
Then  the system
$$
a_iv_i+v_{i+1}=1\quad i=1,\dots,n
$$
has determinant $\prod_{\ell=1}^n a_{\ell}+(-1)^{n-1}$
and a unique solution  given by
$$
v_i=\frac{W_i}{D}:=
\frac{\sum_{j=1}^n(-1)^{j-1}\prod_{\ell=i+j}^{i+n-1}a_{\ell}}
{\prod_{\ell=1}^n a_{\ell}+(-1)^{n-1}}. \qed
$$
\end{lem}

Next we consider in greater detail the 
two low dimensional cases.

\begin{say}[Quasi--smooth rational curves]\label{qs.ratcurves} 

The case $n=3$ gives quasi--smooth rational curves
in weighted projective planes.
Here we have a system
$$
a_1w_1+w_2=a_2w_2+w_3=a_3w_3+w_1=d
$$
with solutions
$$
w_1=\frac{a_2a_3-a_3+1}{w^*},\
w_2=\frac{a_3a_1-a_1+1}{w^*},\
w_3=\frac{a_1a_2-a_2+1}{w^*},\
d=\frac{a_1a_2a_3+1}{w^*}.
$$
We can also compute the genus of  
the general member of the linear system
$$
C\in |x_1^{a_1}x_2, x_2^{a_2}x_3, x_3^{a_3}x_1|
$$
using the adjunction formula (\ref{all.qs.rcs}.4).
By explicit computation
$$
\frac{d(d-w_1-w_2-w_3)}{w_1w_2w_3}=w^*-
\frac1{w_1}-\frac1{w_2}-\frac1{w_3},
$$
Thus  by (\ref{all.qs.rcs}.4),
 the genus of $C$ is  $(w^*-1)/2$.
Thus $C$ is a smooth rational curve iff $w^*=1$,
that is, when 
$a_2a_3-a_3+1, a_3a_1-a_1+1$ and $ a_1a_2-a_2+1$
 are relatively prime.

The Kodaira dimension of the pair $(P:=\p(w_1,w_2,w_3),C)$
is determined by the sign of 
$$
\begin{array}{rcl}
\deg (C+K_{P})&=&a_1a_2a_3+1-(a_1a_2-a_2+1+ a_2a_3-a_3+1+a_3a_1-a_1+1)\\
&=&
(a_1-1)(a_2-1)(a_3-1)-1.
\end{array}
$$
This is negative if one of $a_1,a_2,a_3$ is 1.
If say $a_3=1$ then we get weighted projective planes
$\p^2(a_2,1, a_1a_2-a_2+1)$ with only 2 singular points
and the corresponding link is a lens space.

If $a_1=a_2=a_3=2$ then the relatively prime conditions is not satisfied.
In all other cases $(a_1-1)(a_2-1)(a_3-1)-1>0$ so the Kodaira dimension is
2. \qed
\medskip

As a side remark we note that 
(\ref{H.cob.exmps}) lists all
interesting quasi--smooth rational curves in
weighted projective planes.

\begin{prop}\label{all.qs.rcs}
 Let $P:=\p(u,v,w)$ be a well formed weighted
projective plane and $C=C_d\subset P$ a quasi--smooth
rational curve of degree $d$.
Then, up to permuting the coordinates and
isomorphism,  the pair $d,C$ is one of the following
\begin{enumerate}
\item $d=w$ and $C=(z=0)$.
\item $d=u+v, w|u+v$ and $C\in |xy, z^{(u+v)/w}, \dots|$ where the
existence of other degree $d$ monomials depends on further numerical
coincidences.
\item $(P,C)$ is as in (\ref{H.cob.exmps}).
\end{enumerate}
\end{prop}

Proof. Let $C\subset S$ be a quasi--smooth curve on a surface $S$
which passes through the singular points $P_i$
which are cyclic quotients by $\z/m_i$. The adjunction
formula (cf.\ \cite{cort}) says that
$$
C(C+K_S)=2g(C)-2+\sum \bigl(1-\tfrac1{m_i}\bigr).
\eqno{(\ref{all.qs.rcs}.4)}
$$
Assume now that $S=P$ is a weighted
projective plane, $C$ is rational of degree $d$ and it
passes through at most 2 singular points of indices $u,v$,
where $u=1$ or $v=1$ are allowed. Then we get that
$$
\tfrac{d(d-u-v-w)}{uvw}=-\tfrac1{u}-\tfrac1{v}.
$$
Thus $d<u+v+w$ and so if $z$ is the coordinate with the
biggest weight then it  appears in one of the monomials $z, z^2, zy,zx$.
The rest follows by an easy case analysis, giving the
first two possibilities.

It remains to consider the case when $C$ passes through
all 3 singular points and $u,v,w\geq 2$. This gives the equation
$$
\tfrac{d(d-u-v-w)}{uvw}=1-\tfrac1{u}-\tfrac1{v}-\tfrac1{w}.
$$
Aside from the case $\{u,v,w\}=\{2,3,5\}$, the right hand side
is positive and there is no easy upper bound for $d$.

The quasi--smoothness conditions 
show that, up to permuting the coordinates,
we have monomials
$$
x^ay,y^bz,z^cx\qtq{or}  x^ay, y^bx, z^cx.
$$
The first of these leads to (\ref{H.cob.exmps}) and to our
last possibility. (In fact one can check that
in the 3 singular point case, there are no other monomials
with the same weighted degree.)

Finally we  exclude
the case $x^ay, y^bx, z^cx$. All of these are divisible
by $x$, thus we also must have another monomial
$z^{c'}y$. Thus $(c-c')w=v-u$ and so $w|u-v$.
Up to interchanging $x,y$ we can assume that $u\geq w+1$.

From $d=ua+v=bv+u$ we obtain that $(a-1)u=(b-1)v$. Since
$(u,v)=1$ we get that $d=muv+u+v$ for some $m\geq 1$.
Substituting into the genus formula and rearranging we get
$(mu+1)(mv+1)=(m+1)w$.  But
$mu+1\geq u+1>w$ and $mv+1\geq 2m+1$ give a contradiction.\qed
\end{say}

\begin{say}[Examples of rational homology $\c\p^2$-s]\label{rhcp2.say}

For $n=4$ we get
a surface 
$$
S=S(a_1,a_2,a_3,a_4):=(x_1^{a_1}x_2+x_2^{a_2}x_3+x_3^{a_3}x_4+x_4^{a_4}x_1=0)
\subset \p(w_1,w_2,w_3,w_4),
$$
where the $a_i$ and $w_i$ satisfy a system of equations
$$
a_1w_1+w_2=a_2w_2+w_3=a_3w_3+w_4=a_4w_4+w_1=d
$$
with solutions
$$
w_1=\frac{a_2a_3a_4-a_3a_4+a_4-1}{w^*}, \dots\qtq{and} 
d=\frac{a_1a_2a_3a_4-1}{w^*}.
$$
If $w^*=1$ then
$S(a_1,a_2,a_3,a_4)$ is a rational surface
with 4 quotient  singularities at the coordinate vertices
and with $H_2(S,\q)\cong \q^3$.

Note that $S$ contains the two rational curves
$$
C_1:=(x_1=x_3=0)\qtq{and} C_2:=(x_2=x_4=0).
$$
Both of these are quasi--smooth in $S$.
Thus by the adjunction formula (\ref{all.qs.rcs}.4),
$$
(K_S+C_1)\cdot C_1=-\frac1{w_2}-\frac1{w_4}\qtq{and}
(K_S+C_2)\cdot C_2=-\frac1{w_1}-\frac1{w_3}.
$$
This implies that both curves 
have negative intersection with 
$K_S+(1-\epsilon)(C_1+C_2)$ for $0<\epsilon\ll 1$, and so they are
are extremal rays for
the $K_S+(1-\epsilon)(C_1+C_2)$ minimal model program.
(See \cite{km-book} for an introduction.)
Thus 
$C_1$ and $C_2$ are both contractible
to quotient singularities and we get
rational surfaces
$$
\pi:S\to S^*=S^*(a_1,a_2,a_3,a_4).
$$

If the $\{w_i\}$ are pairwise relatively prime, then
the canonical class of $S$ is
$$
K_S=\o_{\p}(\textstyle{\prod} a_i-1-\textstyle{\sum} w_i).
$$
If the pairwise relatively prime assumption fails
then the general adjunction formula  \cite{cort} says that
$$
K_S+\Bigl(1-\frac1{\gcd(w_1,w_3)}\Bigr)C_1
+\Bigl(1-\frac1{\gcd(w_2,w_4)}\Bigr)C_2
=\o_{\p}(\textstyle{\prod} a_i-1-\textstyle{\sum} w_i).
$$

Note that if $a_1,a_2,a_3,a_4\geq 4$ then
$K_{S^*}$ is ample. 
One can write down an explicit formula for the
self intersection of $K_{S^*}$, but it
is rather complicated. In any case, one sees that
it also converges to 1 as $a_1,a_2,a_3,a_4\to\infty$.
\end{say}

\begin{ack}  I thank V.\ Alexeev, D.\ Auroux, J.\ Amoros, 
W.\ Chen, I.\ Dolgachev,
C.\ Gordon, 
J.\ Keum, T.\ Mrowka, A.\ N\'emethi and D.-Q.\ Zhang
for useful comments  and especially 
M.\ Zaidenberg who called my attention to numerous
related papers.
I am grateful to P.\ Hacking and J.\ Wahl
for allowing me to use their unpublished results.
J.M.\ Johnson helped with computer experiments
that lead to the discovery of various properties of the
examples in Section 5.
Partial financial support was provided by  the NSF under grant number 
DMS-0500198. 
\end{ack}

\bibliography{refs}

\bigskip

\noindent Princeton University, Princeton NJ 08544-1000

\begin{verbatim}kollar@math.princeton.edu\end{verbatim}

\end{document}